\documentclass[11pt]{article}
\usepackage{amsmath}
\usepackage{amssymb}
\usepackage{latexsym}
\usepackage{graphicx}
\usepackage{cite}
\usepackage{mathabx}

\usepackage[utf8]{inputenc}
\usepackage[english]{babel}

\usepackage{amsthm, amstext}
\usepackage{array, amsfonts, mathrsfs}

\theoremstyle{plain}
\newtheorem{thm}{Theorem}

\usepackage{pagecolor}
\definecolor{back}{RGB}{209, 208, 178}

\title{Determination of a wave field in a laterally inhomogeneous medium from boundary data}
\author{
M.~N.~Demchenko\footnote{St.~Petersburg Department of
V.\,A.~Steklov Institute of Mathematics of
the Russian Academy of Sciences, 
27 Fontanka, St.~Petersburg, Russia. E-mail: demchenko@pdmi.ras.ru.\newline
\indent The research was supported by the RFBR grant 17-01-00529-a.}}

\date{}

\begin{document}

\maketitle
\begin{flushright}
{\em\large
\begin{tabular}{r}
Dedicated to Mikhail I. Belishev\\
on the occasion of his jubilee
\end{tabular}
}
\end{flushright}

\medskip

\begin{abstract}
We deal with the Cauchy problem for a 
perturbed wave equation
in the half-plane with data given on a part of the space-time boundary. 
The equation in consideration describes
a wave process in a laterally inhomogeneous medium.
We propose a reconstruction algorithm,
which
is applicable to the problem of determining nonstationary wave field from boundary data arising in
geophysics. 

\smallskip

\end{abstract}

\medskip

\section{Introduction}
Let $u(x,y,t)$ be a solution to the following hyperbolic equation
\begin{align}
  &\partial^2_{t} u - \Delta u + q u = 0,
  \label{wave}
\end{align}
$(x,y)\in {\mathbb R}\times[0,\infty)$ (the half-plane),
$t\in{\mathbb R}$,
and satisfy the conditions
\begin{equation}
  u|_{y=0} = f, \quad \partial_y u|_{y=0} = g.
  \label{data}
\end{equation}
In the present paper, the problem of determination of the solution $u$ from the Cauchy data  $f$, $g$ is investigated.
The following particular case is considered:
the coefficient
$q$ in equation~(\ref{wave}) is assumed to be a function
of $x$. 
Thus equation~(\ref{wave}) is a model of a nonstationary wave process
in a laterally inhomogeneous medium.
We will obtain a relation, which allows determining $u$ at 
any given point $(x_0,y_0,t_0)$ from the values of
the functions $f$, $g$ on a certain bounded set
of points $(x,t)$, 
which depends
on $(x_0,y_0,t_0)$.
We will show that this relation applies in the case when
equation~(\ref{wave}) is satisfied in a cylinder
$\Omega\times{\mathbb R}$, where $\Omega$ is a subset of the half-plane.
The intersection of the boundary $\partial\Omega$
with the boundary of the half-plane is assumed to be
non-empty, and the Cauchy data $f$, $g$ are assumed to be given on some set, whose spatial projection 
lies in this intersection.

The Cauchy problem~(\ref{wave}), (\ref{data}) is ill-posed in the sense of Hadamard~\cite{LRSh}.
However,
the solution $u$ is uniquely determined in some part of the
space-time cylinder dependent on the set, on which 
$f$, $g$ are given.
This follows from the unique continuation property across 
a noncharacteristic surface for equation~(\ref{wave}).
This property is established for various types of linear
partial differential equations~\cite{LRSh, is}.
The most complete results concerning hyperbolic equations
were obtained in~\cite{Tat}.
In particular, these results imply the unique continuation property for equation~(\ref{wave}),
assuming that $q$ is a sufficiently smooth function of
$x$, $y$.

In the case $q\equiv 0$, a number of algorithms for solving
the problem in consideration are known.
One of the pioneering results is that of R.~Courant concerning the
ultrahyperbolic equations in the half-space (see~\cite{K}).
The Cauchy problem for the wave equation in a domain was considered in~\cite{DD18, mn99}.
The inversion formula for the problem in the 
three-dimensional half-space obtained in~\cite{Blag}
allows determining a solution to the wave equation
from the Cauchy data given on a certain
unbounded subset of the space-time boundary.

Problems of determination of solutions to hyperbolic equations from boundary data arise in
geophysics~\cite{Kab}, photoacoustic tomography~\cite{Nat},
tsunami source identification pro\-blems~\cite{vorchev},
and coefficient inverse problems~\cite{bel}.

\section{On the Schr\"odinger operator on the line}
To solve the Cauchy problem~(\ref{wave}), (\ref{data}),
we will deal with the generalized eigenfunction expansion 
of the solution $u$ associated with the
Schr\"odinger operator $\mathcal L = -\partial_x^2 + q(x)$ on the line.
In the case $q\equiv 0$, this expansion is identical to the
Fourier transform in $x$
\[
  \hat u(k, y, t) = \frac{1}{\sqrt{2\pi}}\int_{\mathbb R} e^{-ikx} u(x,y,t)\, dx,
\]
which satisfies the equalities
\[
  \partial_t^2 \hat u - \partial_y^2 \hat u + k^2 \hat u = 0,
  \quad
  \hat u|_{y=0} = \hat f,
  \quad
  \partial_y \hat u|_{y=0} = \hat g.
\]
For arbitrary $q$, we will obtain the same Cauchy problem
with a spectral parameter.
We will solve this problem and apply the inverse transform in $k$, which will result in an expression for $u$
involving the data $f$, $g$.
After that we will show that this expression
involves 
the Cauchy data on a bounded set.

We will need some facts on the Schr\"odinger 
operator~\cite{FDD}.
The latter is well-defined under some restrictions 
on the potential $q$.
We will consider the Schr\"odinger operator only 
with a $C^\infty$-smooth compactly supported real-valued 
potential $q$.
In this case, the operator $\mathcal L$ defined on functions
from $L_2({\mathbb R})$ that have the second derivative
from $L_2({\mathbb R})$ (the Sobolev space $W^2_2({\mathbb R})$)
is a self-adjoint operator in $L_2({\mathbb R})$.
The operator $\mathcal L$ has the absolutely continuous spectrum
of multiplicity two that coincides with the half-line
$\{\lambda\geqslant 0\}$.
Besides, $\mathcal L$ may have the discrete spectrum
consisting of finite number of negative simple
eigenvalues, which will be denoted by
$\{-\varkappa_l^2\}_{l=1}^M$, $\varkappa_l>0$.

Let $A$ be a sufficiently large number such that
${\rm supp}\, q \subset [-A,A]$.
Any solution to Schr\"odinger equation
\begin{equation}
  -\partial_x^2 \varphi + q \varphi = k^2 \varphi \label{schrodinger}
\end{equation}
equals a linear combination of the exponentials
$e^{\pm ikx}$ if $|x|>A$:
\begin{align}
  &\varphi(x) = \alpha_+ e^{ikx} + \alpha_- e^{-ikx}, \quad x > A,
  \label{linexpa}\\
  &\varphi(x) = \beta_+ e^{ikx} + \beta_- e^{-ikx}, \quad x < -A.
  \label{linexpb}
\end{align}
For real nonzero $k$,
we define the function $\varphi_1(x,k)$ 
as the solution to equation~(\ref{schrodinger}) on the whole line that satisfies~(\ref{linexpa}), (\ref{linexpb}) with
$\alpha_-=1/\sqrt{2\pi}$, $\beta_+=0$.
Next we define the function $\varphi_2(x,k)$ as the solution
to equation~(\ref{schrodinger}) that satisfies~(\ref{linexpa}), (\ref{linexpb}) with $\alpha_-=0$, $\beta_+=1/\sqrt{2\pi}$.
For $k\in{\mathbb R}\setminus\{0\}$, the pair of the functions
$\varphi_{1,2}(x,k)$
(as well as $\varphi_{1,2}(x,-k)$) 
forms a basis in the space of the generalized eigenfunctions
corresponding to the point $k^2$ of the absolutely
continuous spectrum.
Further we will consider $\varphi_{1,2}(x,k)$ only for $k>0$.

Along with $\varphi_{1,2}(x,k)$, we introduce the function
$\varphi_0(x,k)$ defined for $x\in{\mathbb R}$ and
$k$ from the finite set of imaginary numbers
$\{i\varkappa_l\}_{l=1}^M$ corresponding to the discrete spectrum.
For $k=i\varkappa_l$, we define
$\varphi_0(\cdot, i\varkappa_l)$ as the normed
eigenfunction of $\mathcal L$ corresponding to
the eigenvalue $-\varkappa_l^2$.
The normed eigenfunction is determined uniquely 
up to a factor of the form
$e^{i\alpha}$, which is chosen arbitrarily.

An arbitrary function $\psi\in L_2({\mathbb R})$
assumes the eigenfunction expansion
$\hat \psi = (\hat \psi_0, \hat \psi_1, \hat \psi_2)$
associated with the operator $\mathcal L$,
\begin{equation}
  \hat \psi_j(k) = \int_{\mathbb R} \psi(x)\, \varphi_j(x,k)^*\, dx,
  \quad j=0,1,2
  \label{trans}
\end{equation}
(here and further $\,*$ means the complex conjugation),
where the functions $\hat \psi_{1,2}(k)$ are defined
for $k>0$, while $\hat \psi_0(k)$ is defined
for $k\in\{i\varkappa_l\}_{l=1}^M$.
This expansion is a unitary transformation acting in the following spaces
\begin{equation}
  L_2({\mathbb R}) \to {\cal H} = {\cal H}_0 \oplus {\cal H}_1 \oplus {\cal H}_2, \quad
  {\cal H}_0 = {\mathbb C}^M,\,\,
  {\cal H}_1 = {\cal H}_2 = L_2({\mathbb R}_+).
  \label{spaces}
\end{equation}
The inverse transformation has the following form
\begin{equation}
  \psi(x) = \sum_{j=1,2} 
  \int_{{\mathbb R}_+} \hat \psi_j(k) \varphi_j(x,k)\, dk + 
  \sum_{l=1}^M \hat \psi_0(i\varkappa_l) \varphi_0(x, i\varkappa_l).
  \label{invtrans}
\end{equation}
The diagonal representation of the operator $\mathcal L$ 
is based on the eigenfunction expansion.
Namely, if $\psi$ belongs to the domain of definition
of the operator $\mathcal L$, and $\tau = \mathcal L \psi$, then
\begin{equation}
  \hat \tau_j(k) = k^2 \hat \psi_j(k), \quad j=0,1,2.
  \label{flg}
\end{equation}
It can be seen from the definitions given above
that in equalities of the form~(\ref{flg}),
the variable $k$ can be any positive number for
$j=1,2$, whereas $k\in\{i\varkappa_l\}_{l=1}^M$ for $j=0$.

For a function $\Phi(\lambda)$ of real variable,
the function of the operator $\Phi(\mathcal L)$
can be described by the following equality
\begin{align}
  \hat \tau_j(k) = \Phi(k^2) \hat \psi_j(k), \quad j=0,1,2,
  \label{Phi_L}
\end{align}
where $\tau = \Phi(\mathcal L) \psi$.
We will also use the following representation for
$\Phi(\mathcal L)$ in the case when $\Phi(\lambda)$ 
is a complex-analytic function
\begin{align}
  \Phi(\mathcal L) = \frac{1}{2\pi i} \int_\Gamma 
  \Phi(\lambda) (\mathcal L-\lambda I)^{-1}\, d\lambda,
  \label{Phi_Lres}
\end{align}
where the integral is taken over an appropriate 
clockwise oriented contour
$\Gamma$ that embraces the spectrum of the operator $\mathcal L$.
The resolvent $(\mathcal L-\lambda I)^{-1}$ in this representation
is an integral operator,
 the Schwartz kernel of which
equals the Green's function $G_\lambda(x_0, x)$
(which is continuous) of equation~(\ref{schrodinger}) for $\lambda=k^2$.
Further we will use the following estimate (see, e.g., \cite{FDD})
\begin{align}
  |G_{k^2}(x_0,x)| \leqslant
  C(q)\, e^{-{\rm Im}\, k\, |x-x_0|}/|k|
  \label{Greenest}
\end{align}
valid for any $k\in{\mathbb C}$ such that ${\rm Im}\, k \geqslant 1+ \max\{\varkappa_l\}_{l=1}^M$.

\section{Determination of the solution $u$}\label{determine}
In this section, we will describe the scheme of solution
of the Cauchy problem in con\-si\-de\-ra\-tion assuming that
$q$ is compactly supported.
We will also assume that the solution $u(x,y,t)$
to equation~(\ref{wave}) (and so, the Cauchy data $f$, $g$)
is compactly supported in $x$. Namely,
for any $T>0$, the restriction of
$u$ to the set $\{0\leqslant y \leqslant T,\,|t|\leqslant T\}$ is compactly supported. 
Note that such solutions do exist since one can take
a solution to an initial boundary value problem
for equation~(\ref{wave}) with compactly supported
initial data (note that the given equation
describes waves propagating with finite velocity).
In sec.~\ref{general}, these restrictions on 
the supports of $q$ and $u$ will be eliminated.

Applying the transformation~(\ref{trans}) 
to equation~(\ref{wave}), we obtain ($j=0,1,2$)
\begin{align*}
  0 =
  \int_{\mathbb R} \varphi_j(x,k)^*\,(\partial_t^2 - \partial_y^2 + \mathcal L) u(x,y,t) \,dx 
  =(\partial_t^2 - \partial_y^2 + k^2) \hat u_j(k,y,t).
\end{align*}
Here we switched the order of differentiation with respect to $y$, $t$ and integration, which is justified since
the function $u$ is smooth, $u(\cdot,y,t)$ is compactly supported, and $\varphi_j$ is bounded in $x$ for any fixed $k$.
We also applied relation~(\ref{flg}).
Treating relations~(\ref{data}) in a similar way,
we arrive at the following Cauchy problem
for $\hat u_j(k,y,t)$, $j=0,1,2$:
\begin{equation}
  \partial_t^2 \hat u_j - \partial_y^2 \hat u_j + k^2 \hat u_j = 0,
  \quad
  \hat u_j|_{y=0} = \hat f_j, \quad
  \partial_y \hat u_j|_{y=0} = \hat g_j.
  \label{Cauchyhat}
\end{equation}
The fundamental solution to this problem has the following form
\begin{equation*}
  \theta(y-|t|)\, r(k,y,t), \quad
  r(k,y,t) = \frac{1}{2} J_0\left(ik\sqrt{y^2-t^2}\right)
\end{equation*}
($\theta$ is the Heaviside function, $J_0$ is the Bessel function of the first kind of order $0$).
This can be verified directly 
or derived from~\cite[chap. V]{K}.
Further we will need estimates of the functions
$r(k,y,t)$ and $\partial_y r(k,y,t)$ for large
(generally, complex) $k$.
To obtain these estimates, we use the following
representation of the Bessel function~\cite{K}
\begin{equation*}
  J_0(\zeta) = \frac{1}{\pi} \int_{-\pi/2}^{\pi/2} e^{i \zeta\, {\sin}s} \, ds,
\end{equation*}
which means that
\[
  r(k,y,t) = \frac{1}{2\pi} \int_{-\pi/2}^{\pi/2} 
  e^{k\sqrt{y^2-t^2}\cdot {\sin}s} \, ds =
  \frac{1}{2\pi} \int_{-\pi/2}^{\pi/2} 
  {\rm ch}\left(k\sqrt{y^2-t^2}\cdot {\sin}s\right) ds,
\]
\[
  \partial_y r(k,y,t) = 
  \frac{k y}{2\pi} \int_{-\pi/2}^{\pi/2} 
  \frac{{\rm sh}\left(k\sqrt{y^2-t^2}\cdot {\sin}s\right)}{\sqrt{y^2-t^2}}
  \,\sin s \, ds
\]
(note that $\partial_y r(k,y,t)$ is regular at $|t|=y$).
These relations imply
the desired 
estimates 
for arbitrary $k\in{\mathbb C}$:
\begin{equation}
  |r(k,y,t)| \leqslant \frac{1}{2}\, e^{|k'| \sqrt{y^2-t^2}},
  \quad
  |\partial_y r(k,y,t)| \leqslant C k^2 y\, e^{|k'| \sqrt{y^2-t^2}},
  \label{rest}
\end{equation}
where $k'={\rm Re}\, k$. 
In the second estimate, we used the inequality
$|\zeta^{-1} {\rm sh} \zeta| \leqslant C e^{|{\rm Re}\, \zeta|}$.
Note that these estimates are far from being optimal.
However, they will suffice for our purposes.

Now write the solution to the problem~(\ref{Cauchyhat})
in terms of the fundamental solution
\begin{align*}
  &\hat u_j(k,y_0,t_0) = \frac{1}{2}\left(\hat f_j(k,t_0+y_0) + 
  \hat f_j(k,t_0-y_0)\right) \notag\\
  &+ \int_{|t-t_0|\leqslant y_0} \left[\partial_{y_0} r(k,y_0,t-t_0)
  \hat f_j(k,t) + r(k, y_0, t-t_0)
  \hat g_j(k,t)\right]\, dt.
\end{align*}
Multiplying both sides by $e^{-h k^2}$, $h>0$,
we obtain
\begin{align*}
  &e^{-h k^2} \hat u_j(k,y_0,t_0) = 
  \frac{1}{2} e^{-h k^2} \left(\hat f_j(k,t_0+y_0) + 
  \hat f_j(k,t_0-y_0)\right) \notag\\
  &+ 
  \int_{|t-t_0|\leqslant y_0} e^{-h k^2} \left[\partial_{y_0}r(k, y_0, t-t_0)
  \hat f_j(k,t) + r(k, y_0, t-t_0)
  \hat g_j(k,t)\right]\, dt. 
\end{align*}
For fixed $y_0$, $t_0$,
the function $\hat u_j(k,y_0,t_0)$ of the variable $k$
belongs to ${\cal H}_j$ (see~(\ref{spaces})),
hence both sides of the last equality
belong to ${\cal H}_j$ as well.
Since the index $j$  takes all possible values 
$0,1,2$ both sides of the equality can be considered
as elements of ${\cal H}$.
In view of estimates~(\ref{rest}),
the regularizing factor $e^{-h k^2}$ forces
the integrand on the right hand side to decay rapidly 
for large $k$.
Therefore the integrand can be considered as an element of
the space ${\cal H}$ that depends continuously on $t$,
whereas the integral itself can be considered as that
of a continuous 
${\cal H}$-valued function of $t$.
Now apply the inverse transformation~(\ref{invtrans})
to both sides of the equality.
This transformation applied to the last term 
on right hand side can be passed under the integral sign.
In view of~(\ref{Phi_L}), we obtain the following equality
in $L_2({\mathbb R})$
\begin{align}
  &e^{-h \mathcal L} u(\cdot, y_0, t_0) = 
  \frac{1}{2} e^{-h \mathcal L}  \left(f(\cdot, t_0+y_0)
  + f(\cdot, t_0-y_0)\right) \notag\\
  &+ \int_{|t-t_0|\leqslant y_0} \left[
    \Phi^D_h(\mathcal L; y_0, t-t_0)\, 
  f(\cdot, t) + \Phi^N_h(\mathcal L; y_0, t-t_0)\, g(\cdot, t)\right]\, dt.
  \label{ureghPhi}
\end{align}
Here
\begin{equation}
  \Phi^D_h(\lambda; y, t) = e^{-h\lambda} \partial_y r\left(\sqrt\lambda, y, t\right),
  \quad
  \Phi^N_h(\lambda; y, t) = e^{-h\lambda} r\left(\sqrt\lambda, y, t\right).
  \label{Phi}
\end{equation}
On both sides of~(\ref{ureghPhi}),
the functions of the operator $\mathcal L$ act on the elements
of $L_2({\mathbb R})$.
Note that the corresponding operators are bounded
in $L_2({\mathbb R})$, since
the functions $e^{-h\lambda}$, $\Phi^D_h(\lambda; y, t)$, $\Phi^N_h(\lambda; y, t)$ of $\lambda$ are bounded on the spectrum of $\mathcal L$.

Now we pass to the limit as $h\to 0$ in equality~(\ref{ureghPhi}).
For $\psi \in L_2({\mathbb R})$, the function $e^{-h \mathcal L} \psi$
converges to $\psi$ in $L_2({\mathbb R})$ as $h\to 0$.
Applying this to $\psi = u(\cdot, y_0, t_0)$ and
$\psi = f(\cdot, t_0\pm y_0)$, we obtain
\begin{align}
  &u(\cdot, y_0, t_0) = 
  \frac{1}{2} \left(f(\cdot, t_0+y_0)
  + f(\cdot, t_0-y_0)\right) \notag\\
  &+ \lim_{h\to 0}\int_{|t-t_0|\leqslant y_0} \left[
    \Phi^D_h(\mathcal L; y_0, t-t_0)\, 
  f(\cdot, t) + \Phi^N_h(\mathcal L; y_0, t-t_0)\, g(\cdot, t)\right]\, dt
  \label{uregPhi}
\end{align}
(the limit is understood in the sense of $L_2({\mathbb R})$).

Relation~(\ref{uregPhi}) allows determining the solution $u$
from the Cauchy data $f$, $g$.
This relation is local in $t$ since 
for fixed $y_0$, $t_0$, the function
$u(\cdot, y_0, t_0)$ is determined from the Cauchy data 
given for $t$ ranging over 
the bounded interval $|t-t_0| \leqslant y_0$.
In sec.~\ref{locality}, we will show that
relation~(\ref{uregPhi}) can be localized both in $x$ 
and $t$
in the sense that it is possible 
to determine $u$ at a fixed point $(x_0,y_0,t_0)$
from values of the functions $f$, $g$
on a bounded set dependent on $(x_0,y_0,t_0)$.
However, 
we first verify that relation~(\ref{uregPhi})
is a pointwise equality.
To do this, we consider once more the terms of the form
$e^{-h \mathcal L} \psi$ in equality~(\ref{ureghPhi}).
Observe that $e^{-h \mathcal L} \psi = \Psi(\cdot, h)$, where
$\Psi$ is the solution to the parabolic equation
\[
  \partial_h \Psi = \Delta \Psi - q \Psi,
\]
in which $h$ plays the role of time,
with the initial data $\Psi|_{h=0} = \psi$.
Therefore for any $h>0$, the function $e^{-h \mathcal L} \psi$
is continuous in ${\mathbb R}$.
Moreover, 
in the context of~(\ref{ureghPhi}),
the initial function $\psi$ is smooth and has compact 
support, which implies that
\[
  (e^{-h \mathcal L} \psi)(x_0) = \Psi(x_0,h) \to \psi(x_0), \quad
  h\to 0
\]
for any $x_0\in{\mathbb R}$ (see, e.g., \cite[chap. IV]{Lad}).
From these observations, it follows that the left hand side
and the first term on the right hand side of~(\ref{ureghPhi})
are continuous functions in ${\mathbb R}$,
having pointwise limits as $h\to 0$.
Hence the same is true for the integral on the right hand side of~(\ref{ureghPhi}). 
Equality~(\ref{uregPhi}) now reads as follows
\begin{align}
  &u(x_0, y_0, t_0) = 
  \frac{1}{2} \left(f(x_0, t_0+y_0)
  + f(x_0, t_0-y_0)\right) \notag\\
  &+ \lim_{h\to 0}\left\{\int_{|t-t_0|\leqslant y_0} 
    \left[\Phi^D_h(\mathcal L; y_0, t-t_0)\, 
  f(\cdot, t) + \Phi^N_h(\mathcal L; y_0, t-t_0)\, g(\cdot, t)\right] dt
\,\bigg|_{x_0}\right\}
\tag{\ref{uregPhi}${}^\prime$}
  \label{uregPhi'}
\end{align}
for any $x_0, t_0\in{\mathbb R}$, $y_0>0$.
The expression in the figure brackets is
the value of the integral (the latter is understood as a function in ${\mathbb R}$) at $x_0$.

\section{Localization of relation~(\ref{uregPhi'})}
\label{locality}
Now we investigate the functions of the operator $\mathcal L$
occurring on the right hand side of~(\ref{uregPhi'}).
According to~(\ref{Phi_Lres}), we have ($|t|\leqslant y_0$)
\begin{align*}
  \Phi^N_h(\mathcal L; y_0, t) = \frac{1}{2\pi i} \int_\Gamma 
  \Phi^N_h(\lambda; y_0, t) (\mathcal L-\lambda I)^{-1}\, d\lambda.
\end{align*}
Make a substitution $k=\sqrt\lambda$ in the integral
assuming that
${\rm Im}\, k > 0$ on the contour of integration.
After a deformation of the contour, we obtain
\begin{align*}
  \Phi^N_h(\mathcal L; y_0, t) = \frac{1}{\pi i} 
  \int_{-\infty+i c/h}^{+\infty+i c/h}
  \Phi^N_h(k^2; y_0, t)\, (\mathcal L-k^2 I)^{-1} k\, dk.
\end{align*}
Here $c, h$ are any positive numbers such that
$c/h > \varkappa_l$, $l=1\ldots M$.
For the Schwartz kernel $K^N_h(x_0,x; y_0, t)$
of the operator $\Phi^N_h(\mathcal L; y_0, t)$, this equality
implies
\begin{align*} 
  K^N_h(x_0, x; y_0, t) = \frac{1}{\pi i} 
  \int_{-\infty+i c/h}^{+\infty+i c/h}
  \Phi^N_h(k^2; y_0, t)\, G_{k^2}(x_0, x) k\, dk
\end{align*}
(recall that $G_{k^2}(x_0, x)$ is the Green's function
of equation~(\ref{schrodinger}) for $\lambda=k^2$).
The function $\Phi^N_h$ decays rapidly for large $k$.
Taking also into account that
the kernel $G_{k^2}(x_0, x)$ is bounded on the contour
of integration, which follows from~(\ref{Greenest}),
and is continuous in $x_0$, $x$,
the kernel $K^N_h$ is continuous in $x_0$, $x$ as well.

Making a substitution in the integral in the previously obtained formula for $K^N_h$ and using the definition~(\ref{Phi}),
we obtain the following relation for $|t|\leqslant y_0$
\begin{align} 
  &K^N_h(x_0, x; y_0, t) = 
  \frac{1}{\pi i} \int_{-\infty+i c}^{+\infty+i c}
  e^{-k^2/h}\, r\left(k/h, y_0, t\right)\, G_{(k/h)^2}(x_0, x) h^{-2} k\, dk.
  \label{intk}
\end{align}
For any fixed $c$ and sufficiently small $h$,
the absolute value of the integrand is estimated 
using~(\ref{Greenest}) and~(\ref{rest}) by the expression
(up to the factor $C(q)$)
\begin{align*}
  &h^{-1} e^{-c\,|x-x_0|/h} e^{(|k'|\, z -{\rm Re}\, k^2)/h} =
  h^{-1} e^{-c\,|x-x_0|/h} e^{(|k'|\, z -k'^2+c^2)/h},
\end{align*}
where $z=\sqrt{y_0^2-t^2}$.
Since $e^{|k'|\,z/h} < e^{k' z/h} + e^{-k' z/h}$,
the resulting expression is majorized by
\begin{align*}
  h^{-1} e^{(-c\,|x-x_0| + z^2/4 + c^2)/h} 
  \sum_\pm e^{-(k'\pm z/2)^2/h}.
\end{align*}
Now taking $c = |x-x_0|/2$ and integrating with respect to 
$k'$, we deduce that the absolute value of the 
integral~(\ref{intk}) does not exceed
\[
  C h^{-1/2} e^{(y_0^2 - t^2 - (x-x_0)^2)/(4h)}.
\]
Now it can be seen that in the case $(x-x_0)^2+t^2 > y_0^2$,
the integral~(\ref{intk}) tends to zero as $h\to 0$.
Thus
\begin{equation}
  \lim_{h\to 0} K^N_h(x_0,x; y_0, t) = 0, 
  \quad
  (x-x_0)^2 + t^2 > y_0^2,
  \label{Klim}
\end{equation}
where the limit is uniform with respect to $x$, $y_0$, $t$
belonging to any compact subset of
\[
  \{ (x,y_0,t)\,|\,\, (x-x_0)^2 + t^2 > y_0^2, 
  \,|t| \leqslant y_0\}.
\]

For the kernel $K^D_h(x_0,x; y_0, t)$ of the operator
$\Phi^D_h(\mathcal L; y_0, t)$, the assertion analogous to~(\ref{Klim}) holds true.

It is convenient to introduce the notation
\[
  K_h = \left(K^D_h, K^N_h\right), \quad
  F = \left(\begin{array}{c}
        f \\ g
  \end{array}
  \right).
\]
In this notation,
relation~(\ref{Klim}) and its counterpart for $K^D_h$
takes the following form 
\begin{equation}
  \lim_{h\to 0} K_h(x_0,x; y_0, t) = 0, 
  \quad
  (x-x_0)^2 + t^2 > y_0^2,
\tag{\ref{Klim}${}^\prime$}
  \label{Klim'}
\end{equation}
where the limit is uniform in the sense specified in~(\ref{Klim}).

Now we write relation~(\ref{uregPhi'}) in terms of $K_h$
\begin{align*}
  &u(x_0, y_0, t_0) = \frac{1}{2} \left(f(x_0, t_0+y_0)
  + f(x_0, t_0-y_0)\right)  \notag\\
  &+ \lim_{h\to 0}
  \int_{|t-t_0|\leqslant y_0} dt \int_{\mathbb R} K_h(x_0,x; y_0,t-t_0) F(x, t)\, dx.
\end{align*}
Using relation~(\ref{Klim'}), in which the limit is uniform
on compact sets, and our assumption that the functions
$u$, $f$, $g$ are compactly supported in $x$,
we can replace the integral with respect to $x$ over ${\mathbb R}$
by that over
the interval $|x-x_0|\leqslant y_0+\varepsilon$, $\varepsilon>0$.
Thus we arrive at the relation
\begin{align}
  &u(x_0, y_0, t_0) = \frac{1}{2} \left(f(x_0, t_0+y_0)
  + f(x_0, t_0-y_0)\right)  \notag\\
  &+ \lim_{h\to 0}
  \int_{|t-t_0|\leqslant y_0} dt \int_{|x-x_0|\leqslant y_0+\varepsilon} K_h(x_0,x; y_0,t-t_0) F(x, t)\, dx,
  \label{ulim}
\end{align}
which
is a localized version of relation~(\ref{uregPhi'}), since $u$ is determined from the Cauchy data on the bounded set.

The derivation of~(\ref{ulim}) from~(\ref{uregPhi'}) relies
on the fact that the kernel $K_h$ 
tends to zero as $h\to 0$, if $(x,t)$ lies outside a certain bounded set. 
It should be noted, however, that if $(x,t)$ belongs to this set, then
$K_h$ generally grows exponentially as $h\to 0$.
In the case $q\equiv 0$, this follows from the 
quite explicit formula for this kernel obtained in~\cite{DD18}.
Due to the growth of the integrand in~(\ref{ulim}),
the limit of the integral generally does not exist, if
$f$, $g$ are arbitrary smooth functions that do not correspond to any Cauchy data.
Therefore if the data are known with some error,
it is necessary to approximate this limit with 
the value of the integral computed for some positive $h$.

\section{Generalization of relation~(\ref{ulim})}\label{general}
In this section, we will show that the assumption that 
equation~(\ref{wave}) is satisfied in the entire half-plane is not necessary for the validity of relation~(\ref{ulim}).
We will consider the case where the solution $u(x,y,t)$
is defined for $(x,y)\in\Omega$, $t\in{\mathbb R}$, where
$\Omega$ is a bounded relatively open subset of the half-plane 
$\{y\geqslant 0\}$.
We will assume that the intersection of the boundary $\partial\Omega$
with the line $\{y=0\}$ is non-empty. The Cauchy data
will be assumed to be given on the set, whose spatial projection lies in this intersection.

\begin{figure}[b!]\centering
  \includegraphics[width=.45\textwidth]{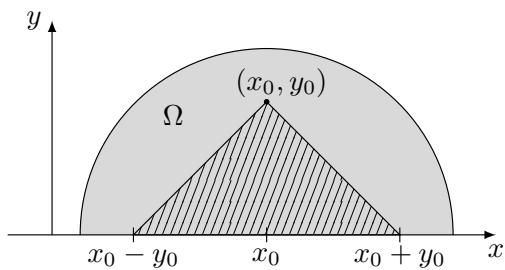}
  \caption{The set $\Omega$ and the set~(\ref{coneset}) (the hatched region).}
  \label{fig}
\end{figure}

Equation~(\ref{wave}) makes sense provided that the coefficient
$q(x)$ is defined on a sufficiently large interval
dependent on $\Omega$.
We will assume that $q(x)$ has a smooth extension to ${\mathbb R}$.
Clearly, this extension can be chosen so that it have compact support.
In this case, the Schr\"odinger operator $\mathcal L$ 
is well defined,
as well as the
kernel $K_h$ occurring on the right hand side of~(\ref{ulim}). 
The following theorem imply, in particular,
that the right hand side in~(\ref{ulim}) does not depend on the choice 
of the extension of the potential $q$. 
\begin{thm}\label{T}
  Suppose that a $C^\infty$-smooth function $u(x,y,t)$ 
  satisfies equation~(\ref{wave}) in the cylinder $\Omega\times{\mathbb R}$.
  The coefficient $q$ in the equation is assumed
  to be a $C^\infty$-smooth compactly supported
  real-valued function of $x\in{\mathbb R}$.
  Let $(x_0,y_0) \in \Omega$ and the closed set
  \begin{equation}
    \{(x,y)\,|\,\, |x-x_0| \leqslant y_0-y, \, y\geqslant 0\}
    \label{coneset}
  \end{equation}
  be a subset of $\Omega$ (see Fig.~\ref{fig}).
  Then for any sufficiently small positive $\varepsilon$,
  relation~(\ref{ulim}) holds true,
  where $f$, $g$ are the Cauchy data determined by
  relations~(\ref{data}).
\end{thm}

To prove the theorem, we choose a smooth 
function $\chi(x,y)$ in the half-plane $\{y\geqslant 0\}$
that has the support in $\Omega$ and is equal to unity
in the neighborhood of the set~(\ref{coneset}).
Put $\tilde u(x,y,t) = \chi(x,y) u(x,y,t)$ if $(x,y)\in\Omega$,
and $\tilde u(x,y,t) = 0$ otherwise.
Thus we obtain a smooth function defined for
$(x,y)$ from the entire
half-plane and for all $t\in{\mathbb R}$.
We have
\begin{align}
  &\partial^2_{t} \tilde u - \Delta \tilde u + q \tilde u = \rho,
  \quad \tilde u|_{y=0} = \tilde f, 
  \quad \partial_y \tilde u|_{y=0} = \tilde g,
  \label{waverti}
\end{align}
where
\[
    \rho = -2\partial_x\chi\, \partial_x u
  - u \Delta\chi, \quad
  \tilde f = \chi f, \quad
  \tilde g = \chi g + f \partial_y\chi
\]
(in the two last equalities, $\chi$ and $\partial_y\chi$ denote
the restrictions of the corresponding functions 
to $\{y=0\}$).

The Cauchy problem~(\ref{waverti}) can be treated 
by the argument 
exposed in sec.~\ref{determine}, \ref{locality}
with the only modification 
concerning 
the presence of $\rho$ on the right hand side of the wave equation. 
This results in the following 
relation 
\begin{align}
  &\tilde u(x_0, y_0, t_0) = \frac{1}{2} \left(\tilde f(x_0, t_0+y_0)
  + \tilde f(x_0, t_0-y_0)\right)  \notag\\
  &+ \lim_{h\to 0} \left\{
  \int_{|t-t_0|\leqslant y_0} dt \int_{|x-x_0|\leqslant y_0+\varepsilon} K_h(x_0,x; y_0,t-t_0) \widetilde F(x, t)\, dx \right.\notag\\
  &\left.+
  \int_0^{y_0} dy \int_{|t-t_0|\leqslant y_0-y} dt 
  \int_{\mathbb R} K^N_h(x_0,x; y_0-y, t-t_0)\, \rho(x, y, t)\, dx\right\},
  \label{ulimti}
\end{align}
which is a substitute for~(\ref{ulim}).
Here $\widetilde F = (\tilde f, \tilde g)^T$, $\varepsilon$ is an arbitrary positive number.
By virtue of relation~(\ref{Klim}),
the integral over ${\mathbb R}$ in the last term in figure brackets can be replaced by
that over any neighborhood of the interval
$|x-x_0| \leqslant y_0-y$.
After this modification, the multiple integral with respect to $x$, $y$, $t$ depends only on the values
of the function $\rho$ at points $(x,y,t)$ such that $(x,y)$
belongs to the corresponding neighborhood of the set~(\ref{coneset}).
However, if this neighborhood is chosen sufficiently small,
then we have $\chi=1$, which implies $\rho=0$.
Thus the term involving $\rho$ in~(\ref{ulimti}) 
can be dropped.
Similarly, $\widetilde F$ can be replaced by $F$
in the first term in figure brackets 
providing that the boundary $\partial\Omega$ contains the interval
$[x_0-y_0-\varepsilon, x_0+y_0+\varepsilon]$
(i.e. $\varepsilon$ is sufficiently small).

It remains to observe that $\tilde f(x_0,t_0 \pm y_0) = f(x_0,t_0 \pm y_0)$ and $\tilde u(x_0,y_0,t_0) = u(x_0,y_0,t_0)$, 
which is due to the equality $\chi=1$ being valid
at the corresponding points.
Thus we arrive at relation~(\ref{ulim}).

The theorem established here means, in particular,
that in the problem on the half-plane considered in sec.~\ref{determine}, \ref{locality}, 
it is possible to eliminate 
the restriction 
that the solution $u$ and the potential $q$ are
compactly supported in $x$.
Indeed, for any given $x_0$, $y_0$,
one can choose a bounded domain $\Omega$ containing the set~(\ref{coneset}) and an arbitrary smooth 
compactly supported potential
 $\tilde q$ such that $\tilde q(x) = q(x)$ whenever 
$(x,y)\in\Omega$.
Then according to Theorem~\ref{T},
the value $u(x_0,y_0,t_0)$ can be found
using relation~(\ref{ulim}), in which
the potential $q$ should be replaced by $\tilde q$.

\end{document}